\documentclass{amsart}

\usepackage{amsthm,amsfonts,amsmath,amssymb,mathrsfs}

\newtheorem{theorem}{Theorem}

\renewcommand{\Re}{\textup{Re}}
\newcommand{\dup}{\textup{d}}
\newcommand{\Int}{\int\limits}
\newcommand{\abs}[1]{\lvert#1\rvert}
\newcommand{\R}{\mathbb R}
\newcommand{\Complex}{\mathbb C}
\newcommand{\la}{\lambda}

\allowdisplaybreaks

\begin{document}

\title{On the generalised Selberg integral of Richards and Zheng}

\author{S. Ole Warnaar}\thanks{Work supported by the Australian Research Council}
\address{Department of Mathematics and Statistics,
The University of Melbourne, VIC 3010, Australia}
\email{warnaar@ms.unimelb.edu.au}



\begin{abstract}
In a recent paper 
Richards and Zheng compute the
determinant of a matrix whose entries are given by beta-type integrals,
thereby generalising an earlier result by Dixon and Varchenko.
They then use their result to obtain a generalisation of the famous
Selberg integral.

In this note we point out that the Selberg-generalisation of
Richards and Zheng is a special case of an
integral over Jack polynomials due to Kadell.
We then show how an integral formula for Jack polynomials of
Okounkov and Olshanski may be applied to prove Kadell's integral
along the lines of Richards and Zheng.
\end{abstract}

\maketitle

Recently, Richards and Zheng established the following theorem \cite{RZ06}.
\begin{theorem}[Richards \& Zheng]\label{thm1}
Let $r$ be a nonnegative integer, $x_1,\dots,x_n\in\R$ and 
$\alpha_1,\dots,\alpha_n\in\Complex$ such that $\Re(\alpha_i)>0$
for all $1\leq i\leq n$.
If 
\[
a_{ij}=\int_{x_i}^{x_{i+1}} y^{j+r-1} 
\prod_{l=1}^n (y-x_l)^{\alpha_l-1} \, \dup y
\]
then
\begin{multline*}
\det_{1\leq i,j\leq n-1}(a_{ij})
=\prod_{1\leq i<j\leq n}(x_j-x_i)
\prod_{\substack{i,j=1 \\ i\neq j}}^n (x_j-x_i)^{\alpha_i-1} \\
\times\sum_{\abs{\nu}=r} \binom{r}{\nu} \:
\frac{\Gamma(\alpha_1+\nu_1)\cdots\Gamma(\alpha_n+\nu_n)}
{\Gamma(\alpha_1+\cdots+\alpha_n+r)} \: x^{r-\nu}.
\end{multline*}
\end{theorem}
In the above we use the following notation.
The sum on the right is over compositions 
$\nu=(\nu_1,\dots,\nu_n)$ of $r$, i.e., $\nu_1,\dots,\nu_n$ are
nonnegative integers such that $\abs{\nu}:=\nu_1+\dots+\nu_n=r$. 
The $\binom{r}{\nu}$ is a multinomial coefficient:
\[
\binom{r}{\nu}=\frac{r!}{(\nu_1)!\cdots(\nu_n)!}\qquad 
\text{for $\abs{\nu}=r$},
\]
and $x^{r-\nu}=x_1^{r-\nu_1}\cdots x_n^{r-\nu_n}$.
Finally we note that the principal branch of each term of the form
$x^{\alpha-1}$ is fixed by $-\pi/2<x<3\pi/2$.

For $r=0$ Theorem~\ref{thm1} is due to Dixon \cite{Dixon05} and
Varchenko \cite{Varchenko89}.

Richards and Zheng use Theorem~\ref{thm1} to prove a generalisation
of the celebrated Selberg integral.
Let $(a)_k$ be a Pochhammer symbol:
\[
(a)_k=a(a+1)\cdots(a+k-1)=\frac{\Gamma(a+k)}{\Gamma(a)},
\]
and define the symmetric polynomial
\[
f_r(x;\gamma)=\sum_{\abs{\nu}=r} \binom{r}{\nu} \:
(\gamma)_{\nu_1}\cdots (\gamma)_{\nu_n}\:
x^{r-\nu}.
\]

\begin{theorem}[Richards \& Zheng]\label{thm2}
For $\alpha,\beta,\gamma\in\Complex$ such that 
\begin{equation*}
\Re(\alpha)>0,~\Re(\beta)>0,~\Re(\gamma)>
-\min\{1/n,\Re(\alpha+r)/(n-1),\Re(\beta)/(n-1)\}
\end{equation*}
and $r$ a nonnegative integer,
\begin{multline*}
\Int_{[0,1]^n}
f_r(x;\gamma) \prod_{1\leq i<j\leq n}\abs{x_i-x_j}^{2\gamma}
\prod_{i=1}^n x_i^{\alpha-1} (1-x_i)^{\beta-1}\, \dup x  \\
=\frac{\Gamma(\alpha) \Gamma(\alpha+\beta+(n-1)\gamma+r)}
{\Gamma(\alpha+r) \Gamma(\alpha+\beta+(n-1)\gamma)} \\ \times
\prod_{i=1}^n
\frac{\Gamma(\alpha+(i-1)\gamma+r) \Gamma(\beta+(i-1)\gamma)\Gamma(i\gamma+1)}
{\Gamma(\alpha+\beta+(i+n-2)\gamma+r)\Gamma(\gamma+1)}.
\end{multline*}
\end{theorem}
Since $f_0(x;\gamma)=1$ this yields the Selberg integral \cite{Selberg44}
when $r=0$. 

Richards and Zheng seem to have been unaware that the polynomial
$f_r(x;\gamma)$ is in fact a Jack polynomial.
Let $g_r^{(\alpha)}(x)$ 
be the symmetric polynomial defined on page 378
of Macdonald's monograph on symmetric functions \cite{Macdonald95}:
\[
\sum_{r=0}^{\infty} g_r^{(\alpha)}(x)\, t^r =
\prod_{i=1}^n (1-tx_i)^{-1/\alpha}.
\]
Using the binomial theorem and the notation $x^{-1}=(x_1^{-1},\dots,x_n^{-1})$
we also have
\[
\sum_{r=0}^{\infty} \frac{f_r(x^{-1};\gamma) (tx_1\cdots x_n)^r}{r!} =
\prod_{i=1}^n\biggl(\;\sum_{\nu_i=0}^{\infty}
\frac{(\gamma)_{\nu_i} (tx_i)^{\nu_i}}{(\nu_i)!}\biggr)
=\prod_{i=1}^n (1-tx_i)^{-\gamma}.
\]
Comparing the above two results we are led to conclude that
\[
f_r(x;\gamma)=r!\, (x_1\cdots x_n)^r g_r^{(1/\gamma)}(x^{-1}).
\]
Now let $P_{\la}^{(\alpha)}(x)$ be the Jack polynomial labelled
by the partition $\la$. On page 380 of \cite{Macdonald95} we find that
\[
P_{(r)}^{(\alpha)}(x)=\frac{r!\, g_r^{(\alpha)}(x)}{(1/\alpha)_r}.
\]
Hence
\[
f_r(x;\gamma)=(\gamma)_r (x_1\cdots x_n)^r P_{(r)}^{(1/\gamma)}(x^{-1}).
\]
Finally we use the well-known fact that for $k$ a positive integer
not exceeding $n$ 
\[
P_{(r^k)}^{(\alpha)}(x)=(x_1\cdots x_n)^r
P_{(r^{n-k})}^{(\alpha)}(x^{-1})
\]
to arrive at
\[
f_r(x;\gamma)=(\gamma)_r\, P_{(r^{n-1})}^{(1/\gamma)}(x).
\]
This result shows that Theorem~\ref{thm2} is in fact the
$\la=(r^{n-1})$ case of Kadell's integral \cite{Kadell97}
(see also \cite[pp. 385--386]{Macdonald95}).
\begin{theorem}[Kadell]\label{thm3}
For $\la$ a partition of at most $n$ parts and
$\alpha,\beta,\gamma\in\Complex$ such that 
\begin{equation*}
\Re(\alpha)>-\la_n,~\Re(\beta)>0,~\Re(\gamma)>
-\min\{1/n,\Re(\alpha+\la_i)/(n-i),\Re(\beta)/(n-1)\}
\end{equation*}
there holds
\begin{align*}
&\frac{1}{n!}\Int_{[0,1]^n}P_{\la}^{(1/\gamma)}(x)
\prod_{1\leq i<j\leq n}\abs{x_i-x_j}^{2\gamma}
\prod_{i=1}^n x_i^{\alpha-1}(1-x_i)^{\beta-1} \; \dup x \\
&\qquad=\prod_{1\leq i<j\leq n}\frac{\Gamma((j-i+1)\gamma+\la_i-\la_j)}
{\Gamma((j-i)\gamma+\la_i-\la_j)} \\
&\quad\qquad \times
\prod_{i=1}^n
\frac{\Gamma(\alpha+(n-i)\gamma+\la_i)\Gamma(\beta+(i-1)\gamma)}
{\Gamma(\alpha+\beta+(2n-i-1)\gamma+\la_i)}.
\end{align*}
\end{theorem}

Given that the Richards--Zheng integral is included in Kadell's
integral it is a natural question to ask for a proof of the latter
along the lines of Richards--Zheng. We address this 
below, giving a very short proof of Theorem~\ref{thm3}.
Key is the next theorem, which follows by taking a limit
in the $q$-integral formula for Macdonald polynomials due to Okounkov
\cite{Okounkov98}
(Okounkov attributes this limiting case to Olshanski).
\begin{theorem}[Okounkov--Olshanski]\label{thm4}
Let $x=(x_1,\dots,x_n)$ and $y=(y_1,\dots,y_{n-1})$ and denote
\[
x_1<y_1<x_2<\cdots<x_{n-1}<y_{n-1}<x_n
\]
by $y\prec x$. For $\la$ a partition of at most $n-1$ parts
\begin{multline*}
P_{\la}^{(1/\gamma)}(x)=
\prod_{i=1}^{n-1}\frac{\Gamma(\la_i+(n-i+1)\gamma)} 
{\Gamma(\la_i+(n-i)\gamma)\Gamma(\gamma)} 
\prod_{1\leq i<j\leq n}(x_j-x_i)^{1-2\gamma} \\
\times \Int_{y\prec x}
P_{\la}^{(1/\gamma)}(y)\prod_{1\leq i<j\leq n-1}(y_j-y_i)
\prod_{i=1}^{n-1}\prod_{j=1}^n \abs{y_i-x_j}^{\gamma-1} \dup y.
\end{multline*}
\end{theorem}
As remarked in \cite{Okounkov98}, for $\gamma=1$ this is equivalent
to the standard definition of the Schur function as a ratio of 
determinants.

\begin{proof}[Proof of Theorem~\ref{thm3}]
First observe that we may assume without loss of generality
that $\la$ has at most $n-1$ parts, i.e., $\la_n=0$.
Indeed, if we define
$\mu=(\la_1-\la_n,\dots,\la_{n-1}-\la_n,0)$
and denote Kadell's integral by $I_{\la}(\alpha,\beta,\gamma)$,
then
\[
I_{\la}(\alpha,\beta,\gamma)=I_{\mu}(\alpha+\la_n,\beta,\gamma)
\]
thanks to
\[
P_{\la}^{(1/\gamma)}(x)=(x_1\cdots x_n)^{\la_n} P_{\mu}^{(1/\gamma)}(x).
\]

We now proceed to prove the theorem by induction on $n$.

For $n=1$ (so that $\la=0$ since $\la_n=0$) the integral is nothing but 
the Euler beta integral \cite{Euler1730}
\[
\int_0^1 x^{\alpha-1} (1-x)^{\beta-1}\, \dup x =
\frac{\Gamma(\alpha) \Gamma(\beta)}{\Gamma(\alpha+\beta)}.
\]

Now assume that $n>1$ and use the symmetry of the integrand to 
write Kadell's integral in the form
\begin{align*}
&\Int_{0<x_1<x_2\cdots<x_n<1}P_{\la}^{(1/\gamma)}(x)
\prod_{1\leq i<j\leq n} (x_j-x_i)^{2\gamma}
\prod_{i=1}^n x_i^{\alpha-1}(1-x_i)^{\beta-1} \, \dup x \\
&\qquad=\prod_{1\leq i<j\leq n}\frac{\Gamma((j-i+1)\gamma+\mu_i-\mu_j)}
{\Gamma((j-i)\gamma+\mu_i-\mu_j)} \\
&\quad\qquad \times
\prod_{i=1}^n
\frac{\Gamma(\alpha+(n-i)\gamma+\mu_i)\Gamma(\beta+(i-1)\gamma)}
{\Gamma(\alpha+\beta+(2n-i-1)\gamma+\mu_i)}.
\end{align*}
Exhibiting the $n$-dependence we denote the
integral on the left by $I_{\la;n}(\alpha,\beta,\gamma)$.

Since $\la$ has at most $n-1$ parts we may eliminate the Jack polynomial 
in the integrand using Theorem~\ref{thm4}. Also interchanging the
order of the $x$- and $y$-integrals then gives
\begin{multline*}
I_{\la;n}(\alpha,\beta,\gamma)=
\prod_{i=1}^{n-1}\frac{\Gamma(\la_i+(n-i+1)\gamma)}
{\Gamma(\la_i+(n-i)\gamma)\Gamma(\gamma)} \\
\times
\Int_{0<y_1<\cdots<y_{n-1}<1}
\!\!
P_{\la}^{(1/\gamma)}(y)\,
J_n(y;\alpha,\beta,\gamma) 
\prod_{1\leq i<j\leq n-1} (y_j-y_i) \; \dup y,
\end{multline*}
where 
\begin{multline*}
J_n(y;\alpha,\beta,\gamma)  =
\Int_{0<x_1<y_1<\cdots<y_{n-1}<x_n<1}
\prod_{1\leq i<j\leq n} (x_j-x_i) \\ \times
\prod_{i=1}^n \biggl[ x_i^{\alpha-1}(1-x_i)^{\beta-1}
\prod_{j=1}^{n-1} \abs{x_i-y_j}^{\gamma-1}\biggr] \dup x.
\end{multline*}

As shown in \cite{RZ02}, an immediate consequence of the
$r=0$ instance of Theorem~\ref{thm1} is the integral
\begin{multline*}
\Int_{y\prec x}
\prod_{1\leq i<j\leq n-1}(y_j-y_i)
\prod_{i=1}^{n-1}\prod_{j=1}^n \abs{y_i-x_j}^{\alpha_j-1} \dup y \\
=\frac{\Gamma(\alpha_1)\cdots\Gamma(\alpha_n)}
{\Gamma(\alpha_1+\cdots+\alpha_n)}
\prod_{1\leq i<j\leq n}(x_j-x_i)^{\alpha_i+\alpha_j-1}.
\end{multline*}
(Peter Forrester pointed out to me that this is also implicit
in Anderson's proof of the Selberg integral \cite{Anderson91}.)
For $\alpha_1=\dots=\alpha_n=\gamma$ this is of course the
$\la=0$ case of Theorem~\ref{thm4}.

The integral $J_n(y;\alpha,\beta,\gamma)$ is precisely of the above
form, with $n\to n+1$ and 
\begin{multline*}
(x_1,\dots,x_n;y_1,\dots,y_{n-1};\alpha_1,\dots,\alpha_{n+1}) \\
\to
(0,y_1,\dots,y_{n-1},1;x_1,\dots,x_n;\alpha,
\underbrace{\gamma,\dots,\gamma}_{\text{$n-1$ times}},\beta).
\end{multline*}
Therefore
\begin{multline*}
J_n(y;\alpha,\beta,\gamma)  =
\frac{\Gamma^{n-1}(\gamma)\Gamma(\alpha)\Gamma(\beta)}
{\Gamma(\alpha+\beta+(n-1)\gamma)} \\ \times
\prod_{1\leq i<j\leq n-1} (y_j-y_i)^{2\gamma-1}
\prod_{i=1}^{n-1} y_i^{\alpha+\gamma-1}(1-y_i)^{\beta+\gamma-1}
\end{multline*}
and, consequently,
\begin{multline*}
I_{\la;n}(\alpha,\beta,\gamma)=
\frac{\Gamma(\alpha)\Gamma(\beta)}
{\Gamma(\alpha+\beta+(n-1)\gamma)} 
\prod_{i=1}^{n-1}\frac{\Gamma(\la_i+(n-i+1)\gamma)}
{\Gamma(\la_i+(n-i)\gamma)} \\
\times
\Int_{0<y_1<\cdots<y_{n-1}<1}
\!\!
P_{\la}^{(1/\gamma)}(y)\,
\prod_{1\leq i<j\leq n-1} (y_j-y_i)^{2\gamma}
\prod_{i=1}^{n-1} y_i^{\alpha+\gamma-1}(1-y_i)^{\beta+\gamma-1}
\; \dup y.
\end{multline*}
Since the right-hand side is Kadell's integral with $n$ replaced by $n-1$
we have
\begin{align*}
I_{\la;n}(\alpha,\beta,\gamma)&=
I_{\la;n-1}(\alpha+\gamma,\beta+\gamma,\gamma) \\ &\qquad
\times
\frac{\Gamma(\alpha)\Gamma(\beta)}
{\Gamma(\alpha+\beta+(n-1)\gamma)} 
\prod_{i=1}^{n-1}\frac{\Gamma(\la_i+(n-i+1)\gamma)}
{\Gamma(\la_i+(n-i)\gamma)} \\
&=
\frac{\Gamma(\alpha)\Gamma(\beta)}
{\Gamma(\alpha+\beta+(n-1)\gamma)} 
\prod_{1\leq i<j\leq n-1}\frac{\Gamma((j-i+1)\gamma+\la_i-\la_j)}
{\Gamma((j-i)\gamma+\la_i-\la_j)} \\
&\qquad \times
\prod_{i=1}^{n-1}
\frac{\Gamma(\alpha+(n-i)\gamma+\la_i)\Gamma(\beta+i\gamma)}
{\Gamma(\alpha+\beta+(2n-i-1)\gamma+\la_i)}\:
\frac{\Gamma(\la_i+(n-i+1)\gamma)}
{\Gamma(\la_i+(n-i)\gamma)},
\end{align*}
where the second equality follows from the usual
induction hypothesis.
It is readily checked that the final expression on the right
is in accordance with the $\la_n=0$ case of Theorem~\ref{thm3},
thus completing the proof.
\end{proof}

\bibliographystyle{amsplain}

\begin{thebibliography}{99}

\bibitem{Anderson91}
G. W. Anderson,
A short proof of Selberg's generalized beta formula,
Forum Math. 3 (1991), 415--417.

\bibitem{Dixon05}
A. L. Dixon,
Generalisations of Legendre's formula $KE'-(K-E)K'=\frac{1}{2}\pi$,
Proc. London Math. Soc. (2) 3 (1905), 206--224.

\bibitem{Euler1730}
L. Euler,
De progressionibus transcendentibus seu quarum 
termini generales algebraice dari nequeunt,
Comm. Acad. Sci. Petropolitanae 5 (1730), 36--57.

\bibitem{Kadell97}
K. W. J. Kadell,
The Selberg--Jack symmetric functions,
Adv. Math. 130 (1997), 33--102.

\bibitem{Macdonald95}
I. G. Macdonald,
Symmetric functions and Hall polynomials,
second edition,
(Oxford University Press, New York, 1995).

\bibitem{Okounkov98}
A. Okounkov,
(Shifted) Macdonald polynomials: $q$-integral representation and 
combinatorial formula,
Compositio Math. 112 (1998), 147--182.

\bibitem{RZ02}
D. Richards, Q. Zheng,
Determinants of period matrices, and an application to 
Selberg's multidimensional beta integral,
Adv. Appl. Math. 28 (2002), 602--633.

\bibitem{RZ06}
D. Richards, Q. Zheng,
The determinant of a hypergeometric period matrix and a
generalization of Selberg's integral, to appear in
Adv. Appl. Math. 

\bibitem{Selberg44}
A. Selberg,
Bemerkninger om et multipelt integral,
Norske Mat. Tidsskr. 26 (1944), 71--78.

\bibitem{Varchenko89}
A. Varchenko,
The Euler beta-function, the Vandermonde determinant,
the Legendre equation, and critical values of linear functions
on a configuration of hyperplanes. I.
Izv. Akad. Nauk. SSSR Ser. Mat. 53 (1989), 1206--1235.

\end{thebibliography}

\end{document}